\documentclass{article}
\usepackage{fancyhdr}
\usepackage{amscd}
\usepackage{graphicx}
\usepackage{amsmath}
\usepackage{amssymb}
\usepackage{amsthm}
\usepackage{mathrsfs,url}
\newtheorem{dfn}{Definition}[section]

\newtheorem{lem}[dfn]{Lemma}

\newtheorem{rem}[dfn]{Remark}
\newtheorem{car}[dfn]{Computer Assisted Result}

\numberwithin{equation}{section}

\usepackage[top=1in, bottom=1in, left=1in, right=1in]{geometry}

\begin{document}
\title{Rigorous numerics of blow-up solutions for ODEs\\ with exponential nonlinearity}

\author{Kaname Matsue\thanks{Institute of Mathematics for Industry, Kyushu University, Fukuoka 819-0395, Japan {\tt kmatsue@imi.kyushu-u.ac.jp}} $^{,}$ \footnote{International Institute for Carbon-Neutral Energy Research (WPI-I$^2$CNER), Kyushu University, Fukuoka 819-0395, Japan}
\and Akitoshi Takayasu\thanks{Faculty of Engineering, Information and Systems, University of Tsukuba, 1-1-1 Tennodai, Tsukuba, Ibaraki 305-8573, Japan ({\tt takitoshi@risk.tsukuba.ac.jp})}
}
\maketitle

\begin{abstract}
Our concerns here are blow-up solutions for ODEs with exponential nonlinearity from the viewpoint of dynamical systems and their numerical validations.
As an example, the finite difference discretization of $u_t = u_{xx} + e^{u^m}$ with the homogeneous Dirichlet boundary condition is considered.
Our idea is based on compactification of phase spaces and time-scale desingularization as in previous works.
In the present case, treatment of exponential nonlinearity is the main issue.
Fortunately, under a kind of exponential homogeneity of vector field, we can treat the problem in the same way as polynomial vector fields.
In particular, we can characterize and validate blow-up solutions with their blow-up times for differential equations with such exponential nonlinearity in the similar way to previous works.
A series of technical treatments of exponential nonlinearity in blow-up problems is also shown with concrete validation examples.
\end{abstract}


\section{Introduction}
Our concerns in this paper are treatments of blow-up solutions for the ODE
\begin{align}
\notag
u_1' &= N^2 (-2u_1 + u_2) + \lambda e^{{u_1^m}},\quad u_{N-1}' = N^2 (u_{N-1}-2u_{N-2}) + \lambda e^{{u_{N-1}^m}},\\
\label{main-eq}
u_i' &= N^2 (u_{i-1} -2u_i + u_{i+1}) + \lambda e^{{u_i^m}},\quad (i=2,\cdots, N-2)
\end{align}
for some positive integer $m$, 
where ${}' = \frac{d}{dt}$ and $\lambda > 0$, from the viewpoint of dynamical systems and their numerical validations.
This system is considered as the finite-difference discretization of the following (initial-)boundary value problem:
\begin{equation}
\label{PDE-original}
\begin{cases}
u_t = u_{yy} + \lambda e^{u^m}, & (t,y)\in (0,T)\times (0,1), \\
u(0,y) = u_0(y), & y\in (0,1), \\
u(t,y) = 0 & \text{at }y=0,1,\quad t\in [0,T)
\end{cases}
\end{equation}
with uniformly grid $y_i = i/N$ and $u_i =u_i(t) \approx u(t,y_i)$.
Dynamics with exponential nonlinearity arises in physico-chemical processing models such as solid fuel ignition and thermal runaway in reaction kinetics \cite{D1985, URP1974} for $m=1$.
In particular, {\em the Arrhenius law} in the theory of reaction kinetics associates the exponential nonlinearity in vector fields. 
The treatment of exponential nonlinearity is therefore a crucial problem when we consider dynamics in realistic scientific arguments.
\par
In preceding works, the authors and collaborators \cite{Mat, MHY2016, MT2017, TMSTMO} have proposed numerical validations of blow-up solutions for {\em polynomial} vector fields as well as their theoretical treatments from the viewpoint of dynamical systems.
The {\em quasi-homogeneity} of vector fields in an asymptotic sense and the {\em time-scale desingularization} are essential tools for geometric treatment of blow-up solutions.
The {\em Lyapunov tracing} technique \cite{MHY2016} realizes the validation of upper and lower bounds of rigorous blow-up times.
However, exponential nonlinearity is against the previous mathematical settings.
Nevertheless, a {\em homogeneous} presence of exponential growth can be reduced to a regular case near infinity as in previous studies, which will give an essence for treating exponential nonlinearity concerning with blow-up behavior both mathematically and numerically.
We believe that the present study is a trigger to treat finite-time singularities including blow-up solutions for vector fields other than polynomial (or rational) ones and their numerical validations.
\par
The organization of this paper is the following.
In Section \ref{section-desing}, we derive a vector field associated with (\ref{main-eq}) appropriate with treatments of blow-up solutions.
The overcome of the presence of exponential nonlinearity is presented there.
In Section \ref{section-procedure}, we show a numerical validation algorithm of blow-up solutions.
The concrete validation results with various aspects in the presence of exponential terms are shown in Section \ref{section-result}.
We end this paper with discussions about future directions in Section \ref{section-discussion}.

\begin{rem}
As for the problem of global existence of solutions for (\ref{PDE-original}) with $m=1$, we can see the classical result in, e.g., \cite{F1969, PV1995}, which shows the existence and non-existence of bounded solutions depending on $\lambda$ and geometry of $\Omega$.
Note that our present study focuses on blow-up computations of spatially discretized problem of PDEs only as a part of ODE problems (cf. \cite{MT2017, TMSTMO}), and this does not directly lead to related results to the original PDEs such as (\ref{PDE-original}) because of difficulties of the present argument in infinite dimensional settings.
Nevertheless, we believe that the present study gives partial aspects of solution structures in associated infinite dimensional problems (i.e., PDEs).
\end{rem}

\section{Desingularized vector field at infinity}
\label{section-desing}

Preceding studies \cite{Mat, MT2017, TMSTMO} begin with choosing appropriate {\em compactifications} for theoretical and numerical studies of blow-up solutions, which embed the phase space into compact manifolds or their tangent spaces.
Typically, compactifications are chosen so that the asymptotically dominant terms at infinity are selected appropriately. 
Such operations can be done for asymptotically {\em polynomial} vector fields and a geometric treatment of blow-up solutions is derived (e.g., \cite{Mat}).
However, the present vector field (\ref{main-eq}) contains exponential nonlinearity, and the general treatment of asymptotic behavior of vector fields at infinity is nontrivial.
Nevertheless, we know that the term $e^{u^m}$ is dominant in (\ref{main-eq}) as $u\to +\infty$ and that exponential terms $\{e^{u_i^m}\}_{i=1}^{N-1}$ are homogeneous in the sense that
\begin{equation}
\label{homogeneity-exp}
e^{(ru_i)^m} = (e^{u_i^m})^{r^m}\quad \text{ for all }r\in \mathbb{R}\quad \text{ and }\quad i=1,\cdots, N-1.
\end{equation}
Therefore, the similar treatment to the polynomial case can be applied to the present problem.
\par
According to the previous works, we choose an appropriate compactification of phase space so that we can treat the infinity with specific direction.
In the present situation, introduce the following {\em directional compactification} (e.g., \cite{Mat}) $\{u_i\}_{i=1}^{N-1}\mapsto \{s, \{x_i\}_{i\in \{1,\dots, N-1 \}\setminus \{N/2\}}\}$ given by
\begin{equation}
\label{compactification}
u_{N/2} = s^{-1},\quad u_i = s^{-1}x_i\quad (i=1,\cdots, N-1,\ i\not =N/2).
\end{equation}
In this case, the infinity corresponds to the subspace $\{s=0\}$ in $(x_1,\cdots, x_{\frac{N}{2}-1}, s, x_{\frac{N}{2}+1, \cdots, x_{N-1}})$-coordinate.
Following \cite{Mat}, we shall call the subspace $\{s=0\}$ {\em the horizon} throughout the paper.
The infinity (in $u_{N/2}$-component) then corresponds to $\{s=0\}$.
\begin{rem}
The choice \eqref{compactification} of directional compactification follows from 
the (numerical) blow-up behavior of \eqref{PDE-original}.
Typical numerical example shows that the solution seems to blow up at the center point $y=\frac{1}{2}$, as seen in Figure \ref{fig_rigorous}, p.8 below.
\end{rem}


Let 
\begin{equation*}
h_{k,\alpha;m}(s) :=  s^{-k}e^{-\alpha / s^m}
\end{equation*}
for integer $k$ and nonnegative real number $\alpha$, and
\begin{equation*}
\begin{cases}
\Delta_{i} := N^2(x_{i-1}-2x_i + x_{i+1}) \quad (i=2,\cdots, N-2, i\not = N/2) & \\
\Delta_{N/2} := N^2(x_{N/2-1}-2 + x_{N/2+1}), \quad 
\Delta_{1} := N^2(-2x_1 + x_2), \quad
\Delta_{N-1} := N^2(x_{N-2}-2x_{N-1}). &
\end{cases}
\end{equation*}
Before our concrete arguments, we check basic properties of $h_{k,\alpha;m}$.
\begin{lem}
\label{lem-h}
For any positive integers $k$, $m$ and nonnegative real number $\alpha$, we have the following properties.
\begin{enumerate}
\item The function $h_{k,\alpha;m}(s)$ is $C^1$ for $s > 0$, and $\lim_{s\to 0+} h_{k,\alpha;m}(s) = 0$.
\item Let
\begin{equation*}
\overline{h_{k,\alpha;m}}(s) := \begin{cases}
h_{k,\alpha;m}(s), & s > 0,\\
0, & s\leq 0.
\end{cases}
\end{equation*}
Then, $\overline{h_{k,\alpha;m}}$ is a $C^1$-extension of $h_{k,\alpha;m}$ over $\mathbb{R}$.
\item We have
$\frac{d}{ds}h_{k,\alpha;m}(s) = h_{k+1,\alpha;m}(s) \left(m\alpha s^{-m}-k\right)$.
\item The function $h_{k,\alpha;m}(s)$ is monotonously increasing over $(0,(m\alpha/k)^{1/m})\subset \mathbb{R}$.
\end{enumerate}
\end{lem}
The proof can be done by direct calculations and l'H$\hat{o}$pital's theorem, and we leave it in the supplemental material \cite{MT_Supp}.
\par
The transformed vector field is
\begin{align*}
s' &= - s \Delta_{N/2} - \lambda (h_{2,1;m}(s))^{-1},\\
x_i' &= - x_i \Delta_{N/2} - x_i \lambda (h_{1,1;m}(s))^{-1} + \Delta_i + \lambda (h_{1,x_i^m;m}(s))^{-1}\quad (i\not = N/2).
\end{align*}
Further introducing the following {\em time-scale desingularization}:
\begin{equation}\label{eqn:timescaling}
\frac{d\tau}{dt} = h_{1,1;m}(s)^{-1},
\end{equation}
we have the following vector field in $\tau$-timescale, which turns out to be regular including $\{s=0\}$:
\begin{align}
\notag
\dot s &= - e^{-1/s^m} \Delta_{N/2} - \lambda s \equiv f_{N/2}(s,x),\\
\label{eqn:desing}
\dot x_i &= - x_i h_{1,1;m}(s)\Delta_{N/2} - x_i \lambda + h_{1,1;m}(s) \Delta_{i} + \lambda h_{0,1-x_i^m;m}(s) \equiv f_{i}(s,x)\quad (i\not = N/2)
\end{align}
where $\dot {} = \frac{d}{d\tau}$.
We shall call it the {\em desingularized vector field} of (\ref{main-eq}).

\begin{rem}
We have chosen the time-scale desingularization $\frac{d\tau}{dt} = T(s)$ so that the following requirements are achieved as in preceding studies \cite{MT2017, TMSTMO}:
\begin{itemize}
\item The desingularized vector field in $\tau$-timescale is smooth (at least $C^1$ in the present argument) on $\{s\geq 0\}$;
\item $T(s)$ is positive for $s>0$ (for orbital equivalence among vector fields);
\item Equilibria on the horizon for the desingularized vector field are non-degenerate (in other words, hyperbolic);
\item $T(s)$ is independent of the remaining $x_i$-components (for simplicity);
\end{itemize}
Note that time-scale desingularizations in preceding studies satisfy all such requirements, which can be proved or verified in theoretical studies \cite{Mat}, although the third requirement depends on problems.
Finally, the smoothness of (\ref{eqn:desing}) including $\{s=0\}$ is ensured by the existence of $C^1$-extensions of $h_{k,\alpha;m}(s)$ over $s\in \mathbb{R}$, as shown in Lemma \ref{lem-h}.
\end{rem}
The Jacobian matrix of the vector field (\ref{eqn:desing}) is given by
\begin{align}
\notag
\frac{\partial f_{N/2}}{\partial s} &= -mh_{m+1,1;m}(s) \Delta_{N/2} - \lambda ,\\
\notag
\frac{\partial f_{N/2}}{\partial x_j} &= -(\delta_{j,N/2-1} + \delta_{j,N/2+1})N^2 h_{0,1;m}(s) \quad (j\not = N/2),\\
\notag
\frac{\partial f_i}{\partial s} &= (1-ms^{-m}) h_{2,1;m}(s) (x_i \Delta_{N/2} - \Delta_{i}) - \lambda m(x_i^m-1) h_{m+1,1-x_i^m;m}(s)\quad (i\not = N/2), \\
\notag
\frac{\partial f_i}{\partial x_j} &= - h_{1,1;m}(s) \Delta_{N/2} \delta_{ij} -N^2 x_i h_{1,1;m}(s) (\delta_{j,N/2-1} + \delta_{j,N/2+1}) - \lambda \delta_{ij}\\
\notag
	&\quad + N^2 h_{1,1;m}(s) \{(1-\delta_{i-1, N/2}) \delta_{i-1, j}(1-\delta_{i-1, 0})  -2\delta_{i,j}\\
\label{Jacobi}
	&\quad + (1-\delta_{i+1, N/2}) \delta_{i+1,j} (1-\delta_{i+1, N}) \} - m\lambda x_i^{m-1} \delta_{i,j} h_{m,1-x_i^m;m}(s)\quad (i,j\not = N/2),
\end{align}
where we have used the (formal) convention $x_{N/2} \equiv 1$ and $\delta_{j,k}$ denotes the Kronecker delta.
Note that exponential nonlinearities appear as {\em exponential decay effects on off-diagonal terms}.
Also, the inequality $x_i \leq 1$ is required for $i=1,\cdots, N-1$ for the boundedness of exponential terms $h_{1,1-x_i^m;m}(s)$.
In the same argument as \cite{Mat}, we can prove that divergent solutions of (\ref{main-eq}) correspond to global trajectories of (\ref{eqn:desing}) in $\{s>0\}$ asymptotic to equilibria (or general invariant sets) on the horizon $\{s=0\}$.
If, moreover, we can prove that maximal existence times $t_{\max}$ of calculated solutions shown below are finite, then the corresponding divergent solutions are actually blow-up solutions.
Therefore, the following procedure can realize rigorous computations of blow-up solutions.

\section{Numerical validation procedure}
\label{section-procedure}
Here, we explain how we validate the rigorous blow-up solutions for \eqref{main-eq} as well as their blow-up times.
As mentioned before, our present methodology is essentially same as that in the preceding works \cite{MT2017, TMSTMO}.
It proceeds in the following steps:
\begin{description}
\item[Step 1.] Determine an equilibrium on the horizon $p_\ast = (0, x_\ast)$.
\item[Step 2.] Validate a \emph{Lyapunov function} $L = L(s,x)$ with $L(0, x_\ast) = 0$ by the procedure given in \cite{MHY2016} as well as its \emph{Lyapunov domain} including $\Omega:=\{(s,x)\in \mathbb{R}_{\ge 0}\times\mathbb{R}^{N-2}\mid L(s,x)\le\epsilon\}$.
\item[Step 3.] Validate $(s(\bar{\tau}),x(\bar{\tau}))\in\mathrm{Int}~\Omega\cap \{s>0\}$.
\item[Step 4.] Obtain rigorous inclusion of the blow-up time.
\end{description}

Firstly, in {\bf Step 1}, an equilibrium on the horizon is obtained by some elementary calculations.
Such a point is the origin of $(s,x)$-coordinate, which is an equilibrium of the desingularized vector field (\ref{eqn:desing}).
Moreover, the equilibrium is asymptotically stable because the Jacobian matrix of the desingularized vector field (\ref{Jacobi}) becomes $\mathrm{diag}([-1,-1,\dots,-1])$ as $s\to 0$.
\par
Secondly, in {\bf Step 2}, we validate the Lyapunov function in the neighborhood of the equilibrium $(s_\ast,x_\ast)=(0,0,\dots,0)$.
From \cite{MHY2016}, the Lyapunov function is given by
\begin{equation}
\label{Lyapunov}
L(s,x)=(s,x)^T(s,x):=s^2+\sum_{i\neq N/2}x_i^2
\end{equation}
with a certain compact set $\bar{\Omega}\subset\mathbb{R}_{\ge 0}\times\mathbb{R}^{N-2}$.
Letting $Df(s,x)$ be the Jacobian matrix of the desingularized vector field at $(s,x)$,
such a $\bar{\Omega}$ is validated by checking whether
\[
	A(s,x):=Df(s,x)^T+Df(s,x),\quad (s,x)\in\bar{\Omega}
\]
is strictly negative definite with \emph{interval arithmetic}.
In practical implementation, we set $\epsilon$ as an upper bound of $L$ so that $\Omega_\epsilon \equiv \{L\leq \epsilon\}\subset\bar{\Omega}$.
\par
Thirdly, in {\bf Step 3}, we rigorously integrate the desingularized ODE until $\tau=\bar{\tau}$ so that $(s(\bar{\tau}),x(\bar{\tau}))\in\mathrm{Int}~\Omega_\epsilon\cap\{s>0\}$ holds.

Finally, in {\bf Step 4}, we rigorously compute the maximal existence time from \eqref{eqn:timescaling}.
It follows
\begin{align*}
	t_{\max}&\equiv\int_{0}^{\infty}s^{-1}e^{-1/s^m}d\tau
	=\int_{0}^{\bar{\tau}}s^{-1}e^{-1/s^m}d\tau+\int_{\bar{\tau}}^{\infty}s^{-1}e^{-1/s^m}d\tau\\
	&\equiv \bar{t}+\int_{\bar{\tau}}^{\infty}s^{-1}e^{-1/s^m}d\tau.
\end{align*}
Because we already obtain the trajectory $(s,x)$ until $\tau=\bar{\tau}$, $\bar{t}$ is rigorously computable by \emph{interval arithmetic}.
The second term is estimated by the \emph{Lyapunov tracing} discussed in \cite{MHY2016,TMSTMO}.
The basic idea is that a fundamental property of Lyapunov function follows
\begin{align}\label{eqn:LT}
	\left.\frac{dL}{d\tau}\left(s(\tau),x(\tau)\right)\right|_{\tau=\bar{\tau}}\le -c_{\bar{\Omega}}L(s(\bar{\tau}),x(\bar{\tau})),
\end{align}
which is strictly negative as long as $(s(\bar{\tau}),x(\bar{\tau}))\not =(s_\ast,x_\ast)$.
Here, $c_{\bar{\Omega}}$ is a positive constant involving eigenvalues of $A(s, x)$ whose details are shown in \cite{TMSTMO}.
From $s\le L(s,x)^{1/2}(\le 1)$ and \eqref{eqn:LT}, we have
\begin{align}
\notag
\int_{\bar{\tau}}^{\infty}s^{-1}e^{-1/s^m}d\tau
\notag
&\le \int_{\bar{\tau}}^{\infty}L(s,x)^{-1/2}e^{-1/L(s,x)^{m/2}}d\tau
\le -\int_{L(s(\bar{\tau}),x(\bar{\tau}))}^{0}L^{-1/2}e^{-1/L^{m/2}}\frac{dL}{c_{\bar{\Omega}}L^{(m+1)/2}}\\
\notag
&=\frac{2}{c_{\bar{\Omega}}m}\int_0^{L(s(\bar{\tau}),x(\bar{\tau}))} \frac{m}{2}L^{-(m/2)-1}e^{-1/L^{m/2}}dL
\notag
=\frac{2}{c_{\bar{\Omega}}m}e^{-1/L(s(\bar{\tau}),x(\bar{\tau}))^{m/2}}\\
\label{eqn:final_extimate}
&\le\frac{2}{c_{\bar{\Omega}}m}e^{-1/\epsilon^{m/2}}.
\end{align}
Now $\bar t$ also has a non-trivial error bound $[\bar{t}_{{\rm low}}, \bar{t}_{{\rm up}}]$ so that $\bar{t}\in [\bar{t}_{{\rm low}}, \bar{t}_{{\rm up}}]$, which is mainly due to rigorous enclosure estimates and rounding errors for validating $\bar{t}$.
The computable bound of $t_{\max}$ is thus given as follows:
\begin{align}
\label{eqn:blowuptime}
	t_{\max}\in \left[\bar{t}_{{\rm low}}, \bar{t}_{{\rm up}}+\frac{2}{c_{\bar{\Omega}}m}e^{-1/\epsilon^{m/2}}\right].
\end{align}

\begin{rem}
\label{rem-Lyap}
Our construction methodology of Lyapunov functions (e.g. \cite{MHY2016}) are based on an appropriate choice of coordinates around equilibria (on the horizon).
The choice can be typically realized by computing eigenvectors of the Jacobian matrices at equilibria, which works successfully {\em if all eigenvalues are simple}.
On the other hand, if there is an eigenvalue with non-trivial multiplicity, we cannot apply the above procedure to Lyapunov function validations because candidates of coordinate transformation can be singular with ordinary numerical procedures.
We easily know that we face this difficulty in the present situation. 
Indeed, $-1$ is the only eigenvalue of the Jacobian matrix with multiplicity $N-1$, as seen in (\ref{Jacobi}).
An alternative way for constructing Lyapunov functions is to apply Schur decomposition of squared matrices, which can be applied regardless of multiplicity of eigenvalues.
Details are shown in \cite{MT2017}.
\end{rem}

The above procedure provides the direct proof of blow-ups for validated solutions {\em in a quantitative way}.
If we only aim at proving that the validated solution blows up, namely without any concrete data for $t_{\max}$, the asymptotic study around equilibria on the horizon based on arguments in \cite{Mat} is applied.
Using the asymptotic method, we can also calculate the blow-up rate of blow-up solutions (see \cite{MT_Supp} for details).
The series of studies reveals {\em a qualitative nature} of validated blow-up solutions.

\section{Numerical validation results}
\label{section-result}

In this section, we show numerical results of our numerical validation procedure.
All computations are carried out on \emph{Bash on Ubuntu on Windows} (ver.\ 16.04), Intel(R) Core(TM) i7-6700K CPU @ 4.00 GHz,  using the \emph{kv library} \cite{kv} (ver.\ 0.4.44) to rigorously compute the trajectories of ODEs.
All codes used to produce the results in this section are freely available from \cite{bib:codes}.

We consider \eqref{main-eq} with the following initial values:
\begin{equation}
\label{initial-1}
u_i(0)=2.5(1-\cos(2\pi y_i))\quad (i=1,2,\dots,N-1)
\end{equation}
for $m=1$, while 
\begin{equation}
\label{initial-2}
u_i(0)= 1-\cos(2\pi y_i)\quad (i=1,2,\dots,N-1)
\end{equation}
for $m=2$.
The corresponding initial data of the desingularized vector field is obtained by
\begin{align*}
s(0)&=\frac{1}{u_{N/2}(0)},\quad 
x_i(0)=\frac{u_i(0)}{u_{N/2}(0)}~\quad (i=1,\cdots, N-1,\ i\neq N/2).
\end{align*}
Our concerning blow-up solutions are trajectories of desingularized vector fields asymptotic to equilibria on the horizon $\{s=0\}$.
As shown in Section \ref{section-procedure}, such equilibria are located at the origin of $(s,x)$-coordinate.
We validate global trajectories of desingularized vector fields asymptotic to the above equilibria with various $N$.
Validated results are collected in Tables \ref{Tab:Ex1} and \ref{Tab:Ex1-m2}.
\begin{table}[t]
\caption{Validated results with $m=1, \lambda = 1$.}
Numerical validations prove that $x(\bar{\tau})\in \mathrm{Int}~\Omega_\epsilon \cap \{s>0\}$.
Furthermore, the rigorous inclusion of the blow-up time $t_{\max}$ is given by the estimate \eqref{eqn:blowuptime} derived in Section \ref{section-procedure}.
\centering
		\begin{tabular}{ccccc}
			\hline 
			$N$ & $\epsilon$ & $\bar{\tau}$ & $t_{\max}$ & Exec. time\\
			\hline\\[-2mm]
			$6$ & $8.02\times 10^{-4}$ & $7.84$ & $0.01223337668427_{7321}^{9155}$ & 4.62 sec.\\[1mm]
			$8$ & $4.11\times 10^{-4}$ & $8.91$ & $0.01384523095580_{1453}^{4485}$ & 9.19 sec.\\[1mm]
			$16$ & $2.81\times 10^{-4}$ & $9.91$ & $0.016198636686_{697263}^{705484}$ & 55 sec.\\[1mm]
			$32$ & $2.32\times 10^{-4}$ & $10.69$ & $0.0167435816193_{27058}^{4389}$ & 8 min. 33 sec.\\[1mm]
			$64$ & $2.32\times 10^{-4}$ & $11.03$ & $0.016874700587_{695561}^{743364}$ & 206 min.\\[1mm]
			$128$ & $2.32\times 10^{-4}$ & $11.44$ & $0.016907356858_{239547}^{428993}$ & 11975 min.\\[1mm]
			\hline 
		\end{tabular}%
	\label{Tab:Ex1}
\end{table}

\begin{table}[t]
	\caption{Validated results $m=2, \lambda = 1$.}
	\centering
	\begin{tabular}{ccccc}
		\hline 
		$N$ & $\epsilon$ & $\bar{\tau}$ & $t_{\max}$ & Exec. time\\
		\hline\\[-2mm]
		$6$ & $7.74\times 10^{-2}$ & $2.97$ & $0.00802832814043_{64432}^{80097}$ & 4.75 sec.\\[1mm]
		$8$ & $7.74\times 10^{-2}$ & $3.33$ & $0.00957793313008_{01847}^{33888}$ & 9.48 sec.\\[1mm]
		$16$ & $2.53\times 10^{-2}$ & $5.11$ & $0.0143478876488_{81462}^{99567}$ & 1 min. 9 sec.\\[1mm]
		$32$ & $2.53\times 10^{-2}$ & $5.69$ & $0.0158714079811_{16483}^{61379}$ & 10 min. 58 sec.\\[1mm]
		$64$ & $1.48\times 10^{-2}$ & $6.61$ & $0.01623711832_{0909481}^{1036679}$ & 242 min. 44 sec. \\[1mm]
		$128$ & $1.48\times 10^{-2}$ & $7.01$ & $0.016328860059_{168021}^{672886}$ & 19049 min. 2 sec.\\[1mm]
		\hline 
	\end{tabular}%
	\label{Tab:Ex1-m2}
\end{table}

\begin{figure}[htbp]
\centering
\begin{minipage}{0.45\hsize}
\centering
\includegraphics[width=5.5cm]{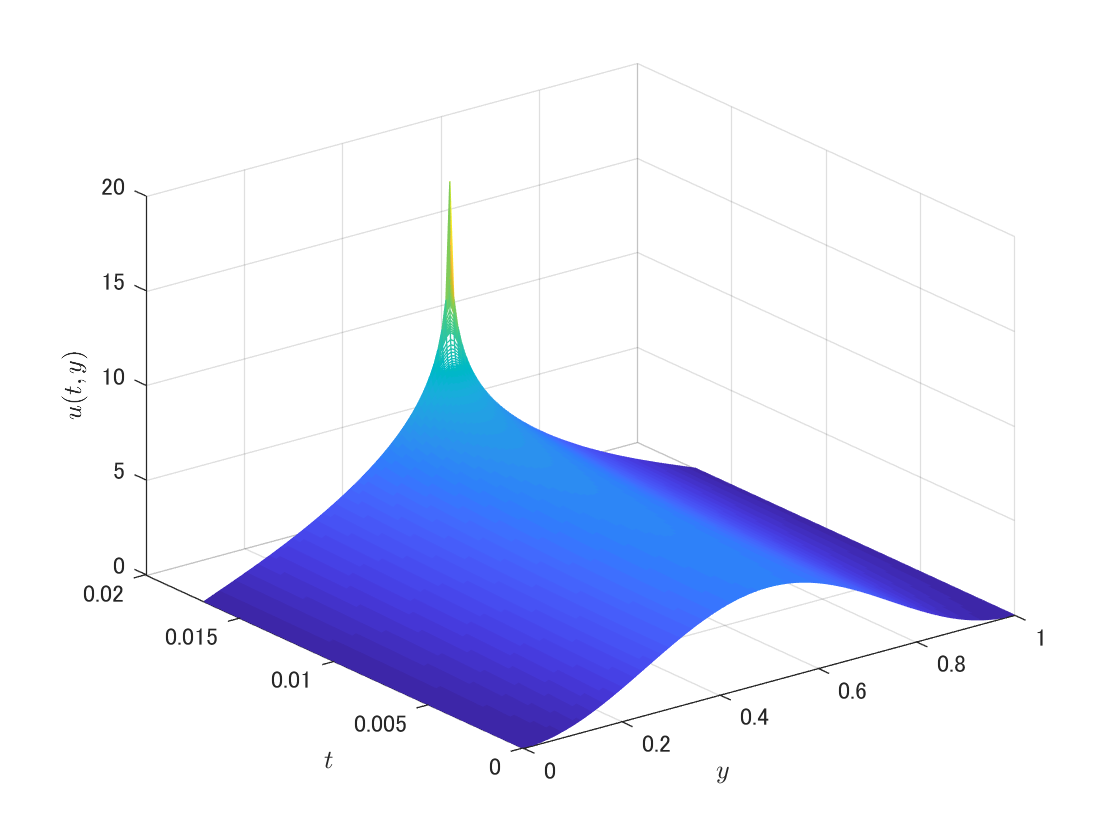}\\
(i)
\end{minipage}
\begin{minipage}{0.45\hsize}
\centering
\includegraphics[width=5.5cm]{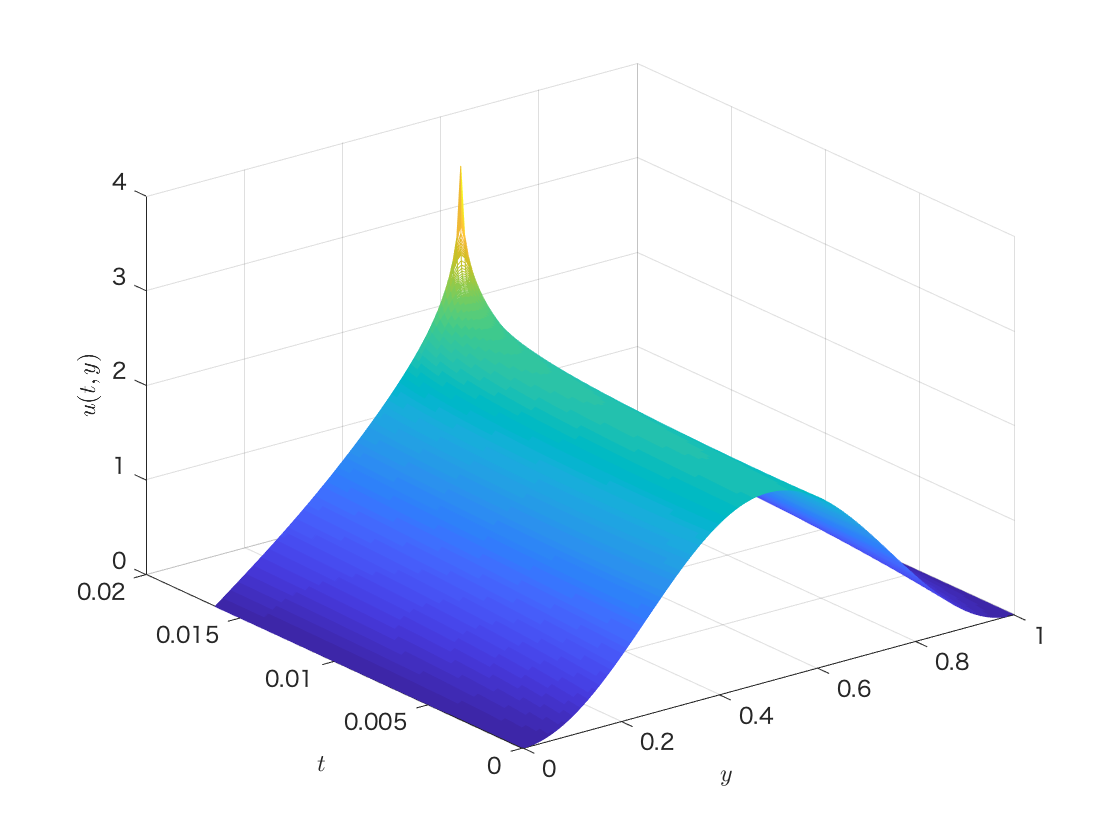}\\
(ii)
\end{minipage}
\caption{Validated blow-up solution profiles for \eqref{main-eq} with $N=128$: $(t,y,u)$-plot. (i) $m=1$. (ii) $m=2$.}
\label{fig_rigorous}
\end{figure}

Applying asymptotic studies of solutions for (\ref{eqn:desing}) 
to the present blow-up solutions, we obtain the following result.
This result characterizes both quantitative and qualitative natures of validated blow-up solutions.
See \cite{MT_Supp} for details of asymptotic behavior of validated solutions.

\begin{car}
Consider (\ref{main-eq}) with $m=1$, $\lambda = 1$, fixed $N$ and initial data (\ref{initial-1}).
Then, for $N=6,8,16,32,64,128$, the solution $\{u_i(t)\}_{i=1}^{N-1}$ blows up at $t = t_{\max} < \infty$, where the $t_{\max}$ is a value inside intervals listed in Table \ref{Tab:Ex1}.
Similarly, replacing $m$ by $2$ and (\ref{initial-1}) by (\ref{initial-2}), we have blow-up validation results with information listed in Table \ref{Tab:Ex1-m2}.
Moreover, all blow-up solutions have the following asymptotic behavior:
$u_{N/2}(t) \sim C\left[\ln\{(t_{\max} - t)^{-1}\}\right]^{1/m}$ as $t\to t_{\max}$ for positive $C$.
\end{car}
Validated blow-up profiles with $N=128$ and $m=1,2$, respectively, are shown in Figure \ref{fig_rigorous}.

\subsection{Practical implementation of enclosures for exponential terms}
In practical computations, we have to care about the treatment of functions of the form $h_{k,\alpha;m}(s) =  s^{-k}e^{-\alpha / s^m}$ near $s=0$.
Actually, the result of numerical computation for $h_{k,\alpha;m}(s)$ becomes infinity  near $s=0$, which is due to {\em zero division}, while the function originally goes to zero as $s\to 0$ for appropriate positive values $(k,\alpha,m)$.
In order to overcome this difficulty, we have applied the monotonous behavior of $h_{k,\alpha;m}(s)$ to implementations for rigorous computations of functions of the form $h_{k,\alpha;m}(s)$ as $s\to 0$.
Observe that, for a positive integers $k, m$ and a real number $\alpha$, $h_{k,\alpha;m}(s)$ is monotonously increasing for sufficiently small but computable $\bar s$ and all $s\in (0,\bar s)$ from Lemma \ref{lem-h}.
Therefore, the function $h_{k,\alpha;m}(s)$ over $[0, \bar{s}]$ is included in the interval $[0, \overline{h_{k,\alpha;m}(\bar{s})}]$ provided $s \in [0,\bar{s}]$, where $\overline{h_{k,\alpha;m}(\bar{s})}$ denotes a computable upper bound of $h_{k,\alpha;m}(\bar{s})$ via interval arithmetic.
As seen below, the exact value of $\overline{h_{k,\alpha;m}(\bar{s})}$ is extremely small for relatively small $\bar{s}$, which helps us with validating targeting objects with very high accuracy.

\subsection{Origin of error bounds for $t_{\max}$}
The value of $e^{-1/\epsilon^{m/2}}$ in the rightmost-hand side of \eqref{eqn:final_extimate} becomes extremely small.
For example, if we put $N = 6$, $m=1$ and $\epsilon = 8.02\times 10^{-4}$, then the upper bound of $t_{\max} - \bar{t}$ is $5.5795557144609417\times 10^{-309}$.
Several typical values of $h_{1,1;m}(\bar s)$ are listed in Table \ref{Tab:Ex2}.
The listed data indicate that $h_{1,1;m}(\bar s)$ is, even for relatively large value $\bar s = 0.02$ in our validations, much less than the machine epsilon $\epsilon_{\rm double} \equiv 2^{-52} \approx 2.22 \times 10^{-16}$ in double-precision arithmetics.
Even in quadruple precision arithmetics, $h_{1,1;m}(\bar s)$ is much less than the corresponding machine epsilon $\epsilon_{\rm quad} \equiv 2^{-112} \approx 1.926 \times 10^{-34}$ for $\bar s = 0.01$.
Therefore, the bound (\ref{eqn:final_extimate}) 
results in multiples of the machine epsilon and the main error of $t_{\max}$ is the difference $\bar t_{\rm up} - \bar t_{\rm low}$; the maximal numerical error of 
\begin{equation*}
\bar{t} = \int_0^{\bar{\tau}} s(\tau)^{-1}e^{-1/s(\tau)^m}d\tau
\end{equation*}
along the component $s(\tau)$ of solution trajectories. 
Note that, in previous studies \cite{MT2017, TMSTMO}, the lower bound of constant $c_{\bar{\Omega}}$ or fractional power of $\epsilon$ corresponding to the rightmost term in (\ref{eqn:final_extimate}) have also caused expansion of error bounds for $t_{\max}$.
In the present case, on the other hand, the exponential decay term drastically decreases such effects.

\begin{table}[t]
	\caption{Values of function $h_{1,1;1}(s)$.}
	\centering
		\begin{tabular}{cc}
			\hline 
			$s$ & $h_{1,1;1}(s)$\\
			\hline\\[-2mm]
			$0.1$ & $4.539992\times 10^{-4}$ \\[1mm]
			$0.05$ & $4.12230724\times 10^{-8}$ \\[1mm]
			$0.02$ & $9.6437492\times 10^{-21}$ \\[1mm]
			$0.01$ & $3.720076\times 10^{-42}$ \\[1mm]
			\hline 
		\end{tabular}%
	\label{Tab:Ex2}
\end{table}

\begin{rem}
Here, we discuss the potential for validating solutions depending on the form of original problems.
In the previous studies, finite-difference discretization of PDEs with polynomial nonlinearity are treated as test examples.
In examples shown in \cite{MT2017, TMSTMO}, 
validations are failed around $N \leq 20$. 
On the other hand, we have succeeded validating blow-up solutions until (at least) $N=128$ in the present study.
We think of the reason very briefly.
The function $h_{k,\alpha;m}(s)$ achieves sufficiently small values for relatively small $s$, in which case the vector field (\ref{eqn:desing}) is approximately governed by
\begin{equation}
\label{eqn:desing-approx}
\dot s \approx - s,\quad 
\dot x_i \approx - x_i\quad (i\not = N/2),
\end{equation}
and the errors are at most on the order of machine epsilon.
Therefore, there are little overestimation due to interval arithmetic, i.e. wrapping effects and any other expanding effects of enclosures due to interactions by other variables during integration of vector fields.
This observation reveals another essence of successful validations in  wide range of $N$, which is a good choice of coordinates (or compactifications).
The equation (\ref{eqn:desing-approx}) is almost diagonal and linear.
According to preceding results in \cite{MT2017, TMSTMO}, if we choose other compactifications (such as Poincar\'{e} type), very limited validation results are obtained.
It is also shown in \cite{MT2017} that, if we choose directional-type compactificaitons like (\ref{compactification}), wider range of validation results are obtained.
Typically, directional-type compactifications and appropriate change of coordinates simplify vector fields near equilibria on the horizon.
The vector field (\ref{eqn:desing-approx}) can be approximately the simplest one in the sense that {\em dynamics of all components are decoupled from each other}, at least near equilibria on the horizon, although it is {\em not} the case of studies in \cite{MT2017}.
Summarizing these observations, {\em good choice of coordinates or compactifications} so that dynamics around equilibria on the horizon becomes as simple as possible, and {\em rapid decay effects} of error terms like exponential decay can extend the range of applicability of proposing blow-up validation methodology.
\end{rem}

\section{Conclusion and Discussions}
\label{section-discussion}

In this paper, we have discussed numerical validations of blow-up solutions for ordinary differential equations with exponential nonlinearity.
The present system is dominated by exponential terms near infinity which treatments in preceding studies (e.g. \cite{Mat, MT2017, TMSTMO}) cannot be directly applied.
Nevertheless, the special homogeneity enables us to apply the similar approach to asymptotically homogeneous cases to studying blow-up behavior of solutions.
A fundamental guideline presented in \cite{Mat, TMSTMO} leads to an appropriate choice of desingularization to the present problem and numerical validation of blow-up solutions.
In the present case, exponential nonlinearity is transformed into exponential {\em decays} for transformed vector fields, which enables us to operate accurate validations.
From the viewpoint of numerical validations, computations of the form $s^{-k}e^{-\alpha/s^m}$ have to be dealt with carefully because of the presence of zero division.
We have applied the monotonous behavior of the function to obtaining the rigorous computable enclosure, which validates various objects effectively.
In particular, our present study has revealed that the theoretical error bound $t_{\max}-\bar t$, presented in (\ref{eqn:final_extimate}), is extremely small because of the exponential decay properties of a priori error bounds.
Exponential decay effects also give a potential to increase the number of discretization $N$ such that numerical validation of blow-up solutions is successful with relatively large $\epsilon$.

\section*{Acknowledgement}
KM was partially supported by World Premier International Research Center Initiative (WPI), Ministry of Education, Culture, Sports, Science and Technology (MEXT), Japan, and JSPS Grant-in-Aid for Young Scientists (B) (No.\ 17K14235).
AT was partially supported by JSPS Grant-in-Aid for Early-Career Scientists (No.\ 18K13453).

\bibliographystyle{plain}
\bibliography{blow_up_exp}

\end{document}